\pdfoutput=1
\documentclass[letterpaper,conference,10pt]{IWMRRF}
\usepackage{amsmath,amssymb,amsfonts,amscd,amsbsy,amsxtra,extarrows}
\usepackage{commath} % for \dif
\usepackage[pdftex]{graphicx}
\usepackage{url}       % Written by Donald Arseneau
\usepackage{amsmath}   % From the American Mathematical Society
\usepackage{array}
\begin{document}

\title{Analytic continuation and differential geometry views on slow manifolds and separatrices}

\author{\authorblockN{Dirk Lebiedz, J\"orn Dietrich, Marcus Heitel, Johannes Poppe}
\authorblockA{Institute for Numerical Mathematics, Ulm University, Germany}}

\maketitle

\begin{abstract}
We start from a mechano-chemical analogy considering the time evolution of a homogeneous chemical reaction modeled by a nonlinear dynamical system (ordinary differential equation, ODE) as the movement of a phase space point on the solution manifold such as the movement of a mass point in curved spacetime. Based on our variational problem formulation \cite{Lebiedz2011} for slow invariant manifold (SIM) computation and ideas from general relativity theory we argue for a coordinate free analysis treatment \cite{Heiter2018} and a differential geometry formulation in terms of geodesic flows \cite{Poppe2019}. In particular, we propose analytic continuation of the dynamical system to the complex time domain to reveal deeper structures and allow the application of the rich toolbox of Fourier and complex analysis to the SIM problem. 
\end{abstract}

\section{Analytic continuation}

\noindent
We introduce analytic continuation of smooth autonomous dynamical systems from real to complex-time domain in order to study slow invariant manifolds (SIM) in terms of spectral properties, geometry and topology of holomorphic curves, respectively embedded Riemann surfaces.
Slow attracting manifolds and, closely related to these, separatrices \cite{Heitel2019} play an important role as backbone structures in phase space and allow to distinguish asymptotic behavior of trajectories. We propose that a holomorphic complex-time view might be useful for identification and analytical or numerical evaluation of characteristic mathematical properties of these flow-invariant structures. 

\noindent
We analyze time-holomorphic solution trajectory manifolds of real-analytic ordinary differential equation (ODE) initial value problems after analytic continuation to the complex time domain:
\noindent
\begin{eqnarray}\label{complexdyn}
\dot{x}&=&f(x), x(t) \in \mathbb{R}^n, t \in \mathbb{R}, x(0)=x_0 \xlongrightarrow{\text{anal. contin.}}  \hspace{1.1cm} \nonumber \\
\noindent
	\dot{z} &=& F(z),
	z=x+iy \in \mathbb{C}^n, t=\sigma +i\tau \in \mathbb{C}, z(0)=z_0 \nonumber
\end{eqnarray}
The real part is $\Re{F}=f$ and the requirement that solutions $z(t)$ are holomorphic functions of complex time $t \in \mathbb{C}$ leads via Cauchy-Riemann differential equations to an ODE for the imaginary time derivative. 

\noindent
The solution trajectories can either be seen as Riemann surfaces defined by images of (subsets of) the complex numbers $\mathbb{C}$ embedded into the ODE solution manifold for a fixed initial value via its complex-time flow or as (ramified or unramified) coverings of the solution manifold. Both viewpoints are geometrical respectively topological in nature and share the potential to access powerful techniques from corresponding mathematical fields, in particular those that relate (complex) analysis with differential geometry and topology. 
Since the phase flow in imaginary time direction is related to complex exponential functions \cite{Dietrich2019}, it turns out to be natural to study Fourier transforms of imaginary time trajectories to analyze spectral properties of the dynamical system. For linear systems classical complex analysis reveals that the spectrum of the system Jacobian is reflected by oscillatory modes in imaginary time direction (see Fig. 1,2 and \cite{Dietrich2019}). We demonstrate for some (linear model and nonlinear Davis-Skodje) benchmark example systems with time scale separation that presence respectively absence of high frequency modes in the Fourier spectrum distinguish initial values on or off slow invariant attracting manifolds of given dimension. 

\noindent
Holomorphic curves respectively Riemann surfaces over complex time allow the definition of a symplectic form related to a surface integral within the solution manifold spanned by dynamical system trajectories. This non-degenerate bilinear 2-form can be discussed in the context of our previous studies of variational formulations of the SIM computation problem in the Hamiltonian framework \cite{Lebiedz2016} whose underlying cotangent bundle geometric structure is symplectic.

\noindent
Imaginary respectively complex time has been discussed in various contexts of modern physical theories such as special and general relativity, quantum mechanics and quantum field theory. 
The time dependent Schr\"odinger equation is formally of parabolic heat-equation type in imaginary time. For time-independent potential energy in the quantum Hamiltonian operator (with the consequence of separability of the time-dependent Schr\"odinger equation), time-differentiation $i \frac{d}{dt}$ can be interpreted as an energy operator. 
Wick rotation is a technique bridging statistical and quantum mechanics by identifying the thermal Boltzmann factor $\frac{1}{kT}$ with imaginary time $\frac{2\pi it}{h}$. In special relativity replacing real by imaginary time makes the Minkowski metric Euclidian. 
For the classical mechanical model of a frictionless pendulum in Hamiltonian formulation, the introduction of imaginary time suggests relations between periodic solutions and the theory of elliptic integrals and functions. 
Obviously, in several contexts there is a relation between physical properties and corresponding modeling techniques and the formal introduction of imaginary/complex time. 
In particular, we show that an extended complex analysis view on flows of dynamical systems sheds new light on the problem of slow invariant manifolds. 

\subsection{Analysis of Riemann surfaces}
\noindent
We demonstrate the value of complex time dynamical systems (analytic ODE) by analyzing a diagonalizable 2-D linear dynamical system with two different eigenvalues (see eq. (1) in \cite{Lebiedz2016} with $\gamma=5$). Figs. 1 and 2 visualize projections of the Riemann surfaces in complex phase space $\mathbb{C}^2$ generated by solution trajectories with initial values off (Fig. 1) and on (Fig. 2) the SIM which is here the slow eigenspace. The oscillatory modulation of the surface can be correlated with an active or relaxed fast mode \cite{Dietrich2019}. 

\begin{figure}[t]
 \centering
 \includegraphics[scale=0.4]
 {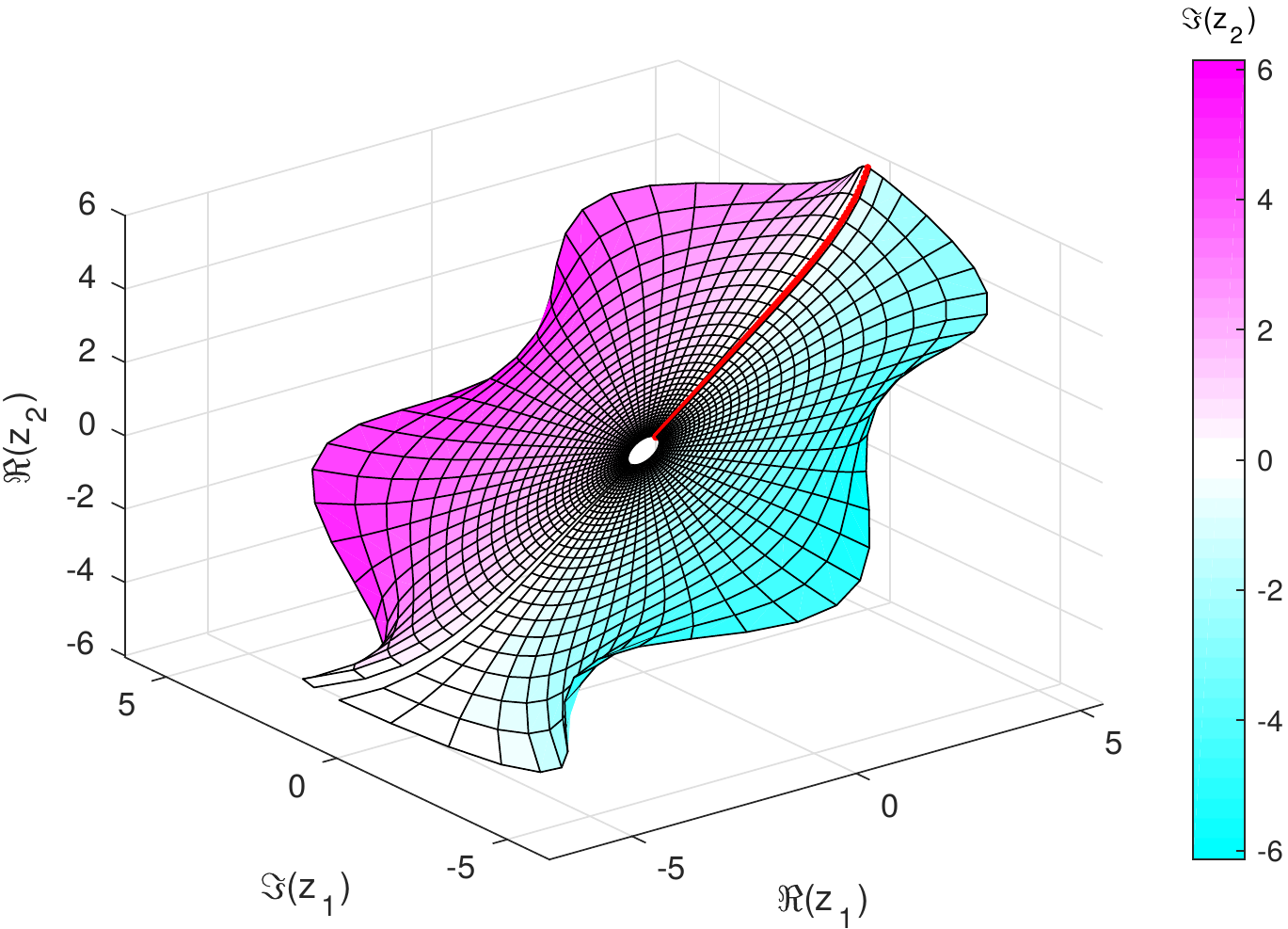}
\caption{2-D linear ODE system (eq. (1) in \cite{Lebiedz2016}, $\gamma=5$) in complex time: Projection of complex solution curve (Riemann surface) for given real initial value $x_0=x(0)$ off the SIM, red line: real time solution trajectory} \label{figure1}
\end{figure}

\begin{figure}[t]
 \centering
 \includegraphics[scale=0.4]
 {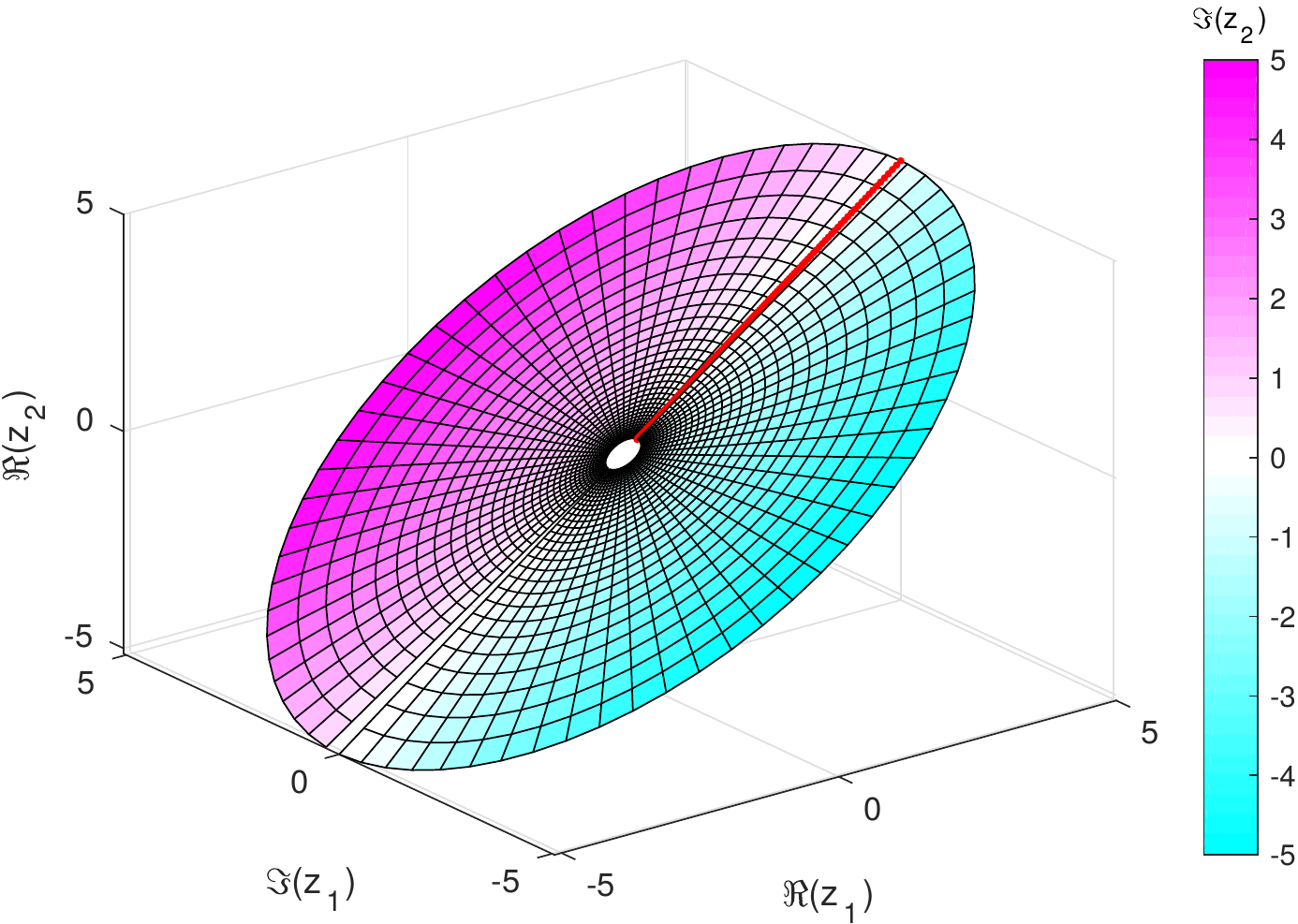}
\caption{2-D linear ODE system (eq. (1) in \cite{Lebiedz2016}, $\gamma=5$) in complex time: Projection of complex solution curve (Riemann surface) for given real initial value $x_0=x(0)$ on the SIM, red line: real time solution trajectory} \label{figure2}
\end{figure}

\subsection{Symmetry considerations}
\noindent
Based on the conception that 1-D SIMs seem to 'balance' contraction rates of trajectories from different phase space directions in \cite{Dietrich2018} an analytic continuation of real to complex phase space is applied to investigate the significance of SIM symmetry properties (Fig. 3) exploiting Noether's theorem. This might help to restrict the potential choice of a suitable Lagrangian function in a variational setting \cite{Lebiedz2016} formulated as an inverse problem for SIM identification. 

\begin{figure}[t]
 \centering
 \includegraphics[scale=0.4]
 {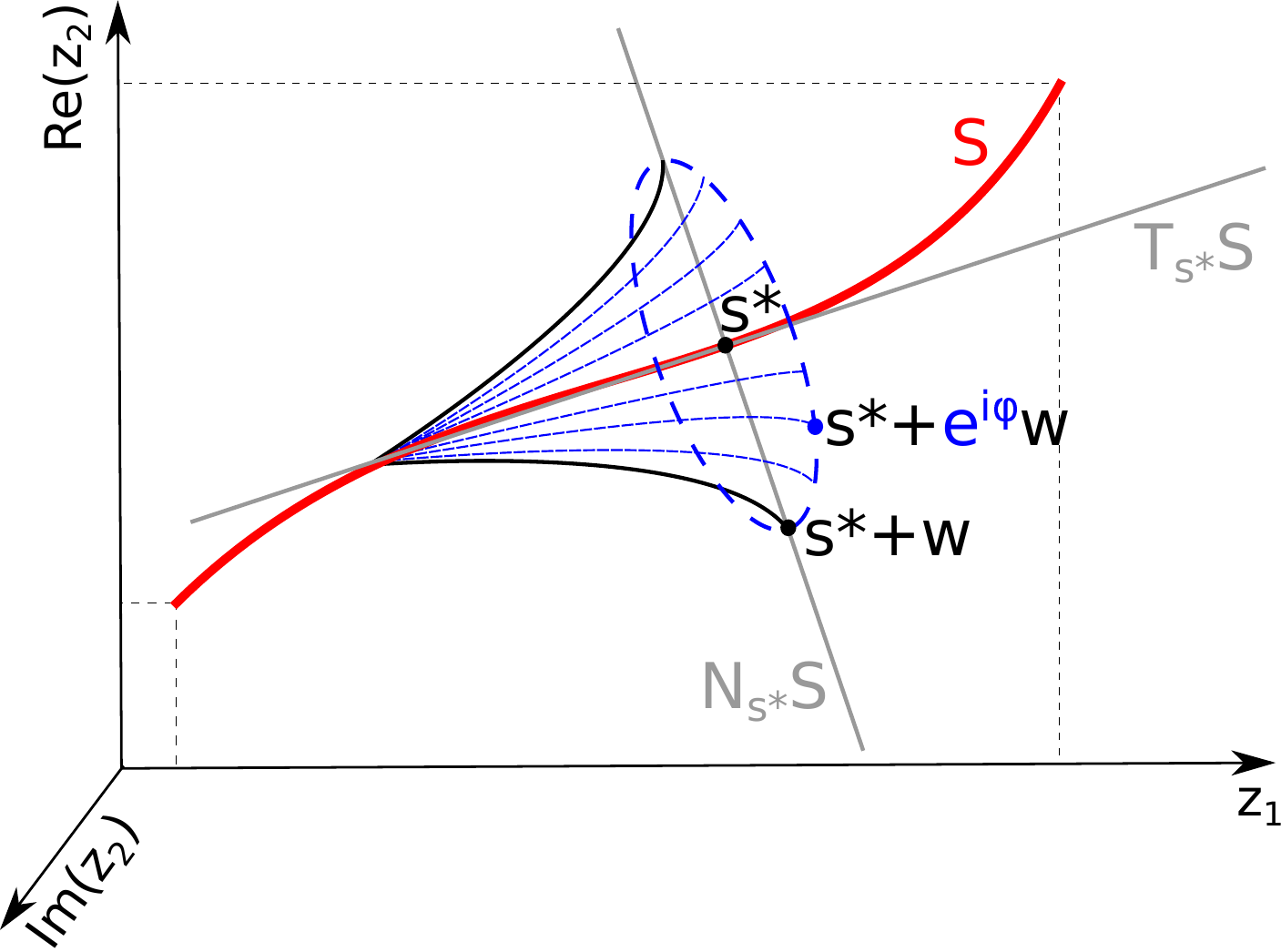}
\caption{Noether's theorem and continuous symmetry (rotation Lie group) of trajectories near SIM in complex phase space (taken from \cite{Dietrich2018})} \label{figure2}
\end{figure}

\section{Differential geometry and geodesic flow}
\noindent
Invariant manifolds in dynamical systems are intrinsic mathematical objects whose characteristics do not depend on a chosen coordinate system (e.g. reaction progress variables for parameterization). 
We start from the ideas to treat the slow manifold problem in the fully coordinate independent setting of differential geometry and the result that necessary (and in special cases also sufficient) conditions for a SIM can be formulated in terms of tensor analysis. The invariance equation, e.g., is reflected by vanishing time-sectional curvature in extended phase space including time as a coordinate axis \cite{Heiter2018}. 
In a second step we derive here a Riemannian metric which makes the solution trajectories geodesic flow lines with respect to the Levi-Civita connection induced by the metric \cite{Poppe2019}. 
In this setting the stretching-based analysis of Adrover et al. \cite{Adrover2007} can be recast in the language of manifold curvature concepts, see \cite{Poppe2019}.

\section*{Acknowledgment}
% optional entry into table of contents (if used)
%\addcontentsline{toc}{section}{Acknowledgment}
\noindent
The Klaus-Tschira foundation is gratefully acknowledged for financial funding of the project.

\end{document}